\documentclass{article}

\def \A {{\mathcal {A}}}

\title{Continuum  pairwise disjoint automorphisms with
Lebesgue spectrum}
\author{ Valery V. Ryzhikov}
\date{}

\begin{document}

\maketitle

\Large
\begin{abstract} A continuous family of pairwise disjoint Gaussian automorphisms is presented, whose even factors have the same Lebesgue spectrum.

\vspace{2mm}
\end{abstract}
\section{Introduction}
Among the automorphisms of a probability space, there is a continuum of pairwise nonisomorphic automorphisms with the same Lebesgue spectrum. This family is formed by $K$-automorphisms with different entropy (a consequence of Kolmogorov's classical results). $K$-automorphisms have no nontrivial factor with zero entropy. Recall that a factor is the restriction of an automorphism to the invariant sigma-algebra of measurable sets. Every factor has a realization as an automorphism as a transformation acting on points. Such a realization can also be called a factor of the original automorphism. For every two $K$-automorphisms, they have isomorphic factors 
 (the necessary information can be found in the book \cite{KSF}).

The strongest form of distinction between ergodic automorphisms $S$ and $S'$ is their disjointness. It means that two factors of an ergodic system, isomorphic to $S$ and $S'$, respectively, are always independent. Pinsker, for example, established that a zero-entropy automorphism and a 
$K$-automorphism are disjoint.

We will consider Gaussian automorphisms of $G(T)$, which are suspensions over certain Sidon transformations of $T$ from \cite{24}.
Thus, we will show that the automorphisms of $G(T)$ are pairwise disjoint, and their even factors $Ev(T)$ have  Lebesgue spectrum. Even factors were considered by Newton and Parry in \cite{NP}. We do not use their definition in our discussion, limiting ourselves to information about the structure of the spectrum of $Ev(T)$. The Gaussian automorphism $G(T)$ as an operator is isomorphic to the direct
sum of operators $T^{\odot n}$, and $Ev(T)$ is isomorphic to the direct
sum of operators $T^{\odot 2n}$. Here $T^{\odot n}$ denotes the symmetric tensor power of order $n$ of the operator $T$. If $\sigma$ is the spectral measure of the operator $T$ (of maximal spectral type), then the convolution power $\sigma^{\ast n}$ of the operator $T^{\odot n}$ is the spectral measure of the operator $T$.
We will consider Sidon transformations of $T$ such that the product $T\times T$ is dissipative, so the spectrum of operators $T^{\odot 2n}$ is Lebesgue.
The spectra of Sidon transformations of $T$ are singular. In this case, the entropy of the Gaussian automorphisms of $G(T)$ is zero.

Thus, in our situation, the even factors of $Ev(T)$ have the same countably multiple Lebesgue spectrum and zero entropy. That's all we need to know about them.

For the proposed Sidon transformations $T$, the automorphisms of $G(T)$ are pairwise disjoint. What metric invariants distinguish them? This role is played by the Kirillov-Kushnirenko entropy. We consider a convenient special version of it, called P-entropy, from the word progession.
It should be noted that A.G. Kushnirenko, thanks to sequence entropy, as far as I know,  was the first to distinguish flows with zero classical entropy and the same Lebesgue spectrum \cite{Ku}.

The paper \cite{RT} gives a generalization of Pinsker's theorem: an automorphism with completely positive P-entropy is disjoint from an automorphism with zero P-entropy. This fact establishes the pairwise disjointness of a specially constructed family of Gaussian automorphisms.

Following the program described above, we obtain a continuous family of pairwise disjoint Newton-Perry factors of Gaussian suspensions over Sidon transformations. These factors have the same Lebesgue spectrum but different Kirillov-Kushnirenko entropies. Our aim is to construct a suitable family of Sidon transformations.

\section{  $P$-entropy of   automorphisms}
We set $P_k=\{n(k),2n(k),\dots, L(k)n(k)\}$, $ L(k)\to\infty.$
For an automorphism $G$ and a finite measurable partition $\xi$ of the probability space, we define the quantities
$$h_k(G,\xi)=\frac 1 {|P_k|} H\left(\bigvee_{p\in P_k}G^p\xi\right),$$
where $H(\xi)$ is the entropy of the partition $\xi$. Recall that for $\xi=\{C_1,C_2,\dots, C_n\}$
$$ H(\xi)=-\sum_{i=1}^n \mu( C_i)\ln \mu( C_i).$$
The $P$-entropy of an automorphism $G$ with respect to $\xi$ is
$$h_{P}(G,\xi)={\limsup_k} \ h_k(G,\xi).$$
And the $P$-entropy of an automorphism $G$ is
$$h_{P}(G)=\sup_\xi h_{P}(G,\xi),$$
where the supremum is taken over all finite measurable partitions.

The metric entropy $h(G)$ of an automorphism $G$ is a special case of $P$-entropy when $P_k=\{1,2,\dots, k\}$.

The equality $h(G^n)=n h(G)$ and the separability of the space of finite partitions obviously imply the following assertion
(see also \cite{RT}).

\vspace{2mm}
\bf Lemma 1. \it If $h(G)=0$, then for every sequence $n(k)\to\infty$ there exists a sequence $ \tilde L(k)\to\infty$ such that for $ L(k)\geq \tilde L(k)$ we have $h_P(G)=0$ for any sequence $P=\{P_k\}$, $P_k=\{n(k),2n(k),\dots, L(k)n(k)\}$.
\rm

\vspace{3mm}

Our goal is to indicate a family of Gaussian automorphisms such that for every pair of its elements $G, G'$, there exists a sequence $P$, as in Lemma 1, such that $h_P(G)=0$, but the automorphism $G'$ has completely positive $h_P$-entropy.
Then, Theorem 4.1 \cite{RT} will imply the disjointness of such $G,G'$.

\section{ Completely positive $P$-entropy of  Gaussian automorphisms}
Let $(R^\infty,\mu)$ be the standard Gaussian probability space.
Orthogonal operators $T$ preserve the measure $\mu$ and are therefore considered automorphisms of $G(T)$ of the probability space, called Gaussian automorphisms (we also call them suspensions).
Gaussian actions have long been studied in ergodic theory \cite{KSF}, 
 and the theory of unitary representations \cite{N}.
An operator $T$ can be considered an orthogonal operator on
an arbitrary infinite-dimensional real Hilbert space isomorphic to the space $R^\infty$ with the standard inner product. For our purposes, we will use operators acting in $L_2(X,m)$, where $(X,m)$ is the standard space with sigma-finite measure $m$. They will be induced by invertible, measure-preserving transformations $T:X\to X$.

\vspace{2mm}
\bf Lemma 2. \it Let the union of increasing finite-dimensional subspaces $H_k$, $k=1,2,\dots,$ of a real Hilbert space $L_2(X,m)$ be dense in it. If there exist sequences $n(k),\, L(k)\to\infty$ such that for an orthogonal operator $T$ the subspaces
$$T^{n(k)}H_k,\,T^{2n(k)}H_k, \dots, T^{n(k)(L(k)-1)}H_k,\,
T^{n(k)L(k)}H_k$$ are pairwise orthogonal for all $k$,
then for $$P=\{P_k\}, \ \ P_k=\{n(k),2n(k),\dots, L(k)n(k)\}$$ the Gaussian suspension $G(T)$ has completely positive $P$-entropy.

\rm

\vspace{3mm}
Proof. Using space isomorphism, we transfer the Gaussian measure from $R^\infty$ to the space $L_2(X,m)$ and denote it by $\mu$. Let $\A_k$ be the algebra of $\mu$-measurable sets whose indicators depend only on vectors in $H_k$.
Then the algebras
$$T^{n(k)}\A_k,\,T^{2n(k)}\A_k,\, T^{n(k)(L(k)-1)}\A_k,\,
T^{n(k)L(k)}\A_k$$ are independent.
This independence obviously implies that the $P$-entropy of the automorphism $G(T)$ is infinite. Moreover,
every $\mu$-measurable finite partition is approximated by $\A_k$-measurable
partitions, so $G(T)$ has completely positive entropy. In fact, the entropy of $G(T)$ with respect to a partition is equal to the entropy of that partition.

Now let us indicate how the operators $T$ that realize completely positive entropy for Gaussian suspensions are structured. Let $X_k$ be an increasing sequence of $m$-measurable sets of finite measure, $\cup_kX_k=X$. Obviously, one can choose a sequence $H_k\subset L_2(X_k,m)$ such that every vector in $L_2(X,m)$ is approximated by vectors in $H_k$ as $k\to\infty$.
If the sets
$$T^{n(k)}X_k,\,T^{2n(k)}X_k, \dots, T^{n(k)(L(k)-1)}X_k,\,
T^{n(k)L(k)}X_k$$
are disjoint, then the spaces $T^{n(k)}H_k,\,T^{2n(k)}H_k, \dots,
T^{n(k)L(k)}H_k$ are pairwise orthogonal. These conditions are easily satisfied, but
we must simultaneously ensure that the spectrum of the transformation $T$ is singular.
For this purpose, we will consider the class of Sidon constructions of rank one \cite{24}.

\section{ Constructions of transformations. Main result}
Let $n_1=1$, $s_j(1)=10n_j$, $j=1,2,\dots,$ and let $s_j(i)$ be natural numbers satisfying the condition $s_j(i)>10s_j(i-1)$, $0<i\leq r_j$,
$r_j=2^j$. The sequence $n_j$ is defined by 
$$ n_{j+1} =n_jr_j +\sum_{i=1}^{r_j}s_j(i).$$
(In the articles, the height of the tower $X_j$ is denoted by $h_j$;
we replaced it with $n_j$, since $h_j$ appeared in the definition of $P$-entropy).

We define the set $X_j$ by induction as the
union of disjoint half-intervals
$E_j, TE_j, T^2E_j,\dots, T^{n_j-1}E_j.$
This set of half-intervals is called the tower of stage $j$; their union is denoted by $X_j$ and is also called a tower.
On these intervals, except for the last one,
the transformation $T$ maps $T^iE_j$ to $T^{i+1}E_j$ by the usual translation.
We represent $E_j$ as the disjoint union of $r_j$ half-intervals
$E_j^1,E_j^2E_j^3,\dots E_j^{r_j}$ of the same length, where $r_j=2^j$.
For each $i=1,2,\dots, r_j$, we consider the so-called column
$$E_j^i, TE_j^i ,T^2 E_j^i,\dots, T^{n_j-1}E_j^i.$$

Let $s_j(1)=10n_j$ and natural numbers $s_j(i)$ be given
satisfying the condition $s_j(i)>10s_j(i-1)$.
To each column $X_{i,j}$, we add $s_j(i)$ half-intervals of the same length as $E_j^i$.
We obtain a set
$$E_j^i, TE_j^i, T^2 E_j^i,\dots, T^{n_j-1}E_j^i, T^{n_j}E_j^i, T^{n_j+1}E_j^i,\dots, T^{n_j+s_j(i)-1}E_j^i$$
(all these sets are disjoint).
For $i<r_j$, we extend the transformation $T$ on the upper floors
of the superstructured columns so that
$T^{n_j+s_j(i)}E_j^i = E_j^{i+1}$.

We denote $E_{j+1}= E^1_j$. All floors of the superposed columns, as can be seen,
form a tower at stage $j+1$,
consisting of the half-intervals
$$E_{j+1}, TE_{j+1}, T^2 E_{j+1},\dots, T^{n(j+1)-1}E_{j+1}.$$

The definition of the transformation $T$ at stage $j$ is preserved at all subsequent stages. As a result, we obtain the space $X=\cup_j X_j$ and an invertible transformation $T:X\to X$, which preserves the standard Lebesgue measure on $X$.

The measure of $X$ is infinite. The transformations defined above form a subclass of Sidon transformations with singular spectrum \cite{24}. The choice of a rapidly growing sequence $r_j=2^j$ ensures that the product $T\times T$ is dissipative (Theorem 3.1 \cite{24}). This means that the products $T^{\odot 2n}$ have a Lebesgue spectrum, and therefore the spectrum of
$Ev(T)$ is Lebesgue.

\vspace{2mm}
\bf Theorem 4.1. \it There exists a continuum of pairwise disjoint Gaussian automorphisms with  (pairwise disjoint) factors of the same  Lebesgue spectrum.\rm

\vspace{2mm}
Proof. We define a continuum family of transformations $T_a$ for all sequences $a=\{a_i\}$, $a_i\in\{0,1\}$ as follows.
Among the stages $j$, we will choose a strictly monotone sequence
$j(k)$ such that
for each $k$, we will specify a set of $2^k$ constructions of transformations of stage $j(k)$, defined on the towers
$X_{j(k)}^a$. Each such tower for all $a$ is the union of $n(k)$
($n(k)$ does not depend on $a$) disjoint half-intervals of equal length.

We continue the further construction of Sidon constructions in a uniform manner and obtain a set $S_k$ of $2^k$ Sidon transformations with singular spectrum. We fix $\pi_k$ -- a finite $\A_k$-measurable partition. By Lemma 1, there exists $L_k$ such that for $P_k=\{n(k),2n(k),\dots, L(k)n(k)\}$, the inequality
$$h_k(G(T_k), \pi_k)< \ 1 / k$$
holds for all $T_k\in S_k$.
We choose $j(k+1)$ such that this inequality holds for all constructions $T$ identical to some $T_k\in S_k$ up to stage $j(k+1)$.
We define a new set $\tilde S_k$ consisting of constructions $\tilde T_k$
identical to $T_k$ up to stage $j(k)$, and at stage $j(k)+1$
let $s_{j(k)+1}> L(k)(n(k)+1)$. Then the sets
$$T^{n(k)}X_{j(k)},\,\ T^{2n(k)}X_{j(k)},\ \dots,\ T^{n(k)(L(k)-1)}X_{j(k)},\, \ T^{n(k)L(k)}X_{j(k)}$$ are disjoint, which leads to
$h_k(G(\tilde T_k), \pi_k)=H(\pi_k).$

Continuing the construction in this way, for every sequence $a$
we obtain the construction $T_a$. If sequences $a$ and $a'$ differ on an infinite set $K=\{k:\, a_k\neq a_k'\}$, then $G(T_a)$ and $G(T_{a'})$ are disjoint. Indeed, let, for example, the set
$K_0= \{k:\, a_k=0, a_k'=1\}$ be infinite. Then, with respect to the sequence of progressions $P=\{P_k\}$, $ k\in K_0$, the automorphism $G(T_{a'})$ has completely positive $h_P$-entropy, but $ h_P(G(T_a))=0$. Thus, the automorphisms $G(T_{a})$ and $G(T_{a'})$
are disjoint. For all $T_{a}$, we ensure that the product $T\times T$ is dissipative. From the set of all $T_a$, we select the desired continuum class and thus complete the proof of the theorem.

\vspace{2mm}
\bf Remark. \rm Instead of a Gaussian suspension $G(T)$, we can consider Poisson suspensions $P(T)$ for a Sidon transformation $T$. The spectra of  $P(T)$ and  $G(T)$ are the same. The entropy of such $P(T)$ is zero, so similar arguments apply to them when constructing
large disjoint families of suspensions $P(T)$ with the same  Lebesgue component in their spectra. However, the situation with factors is  different, $P(T)$ can be prime, see \cite{PR}. It is of interest to  study suspensions for conservative products $T\times T$ with Lebesgue spectrum. This case arises for the Sidon constructions with parameters $r_j$ that grow at a suitable rate, see \cite{26}. It will be nice to have $P(T\otimes T)$
of zero entropy and $G(T\otimes T)$ with infinite one. 

\vspace{4mm}
\ \ \ \ \ \ \ \ \ \ \it $G(S)$ is Bernoulli, $P(S)$ not,  is this  possible?\rm

\large

\end{document}